\documentclass[journal,twoside,web]{ieeecolor}
\usepackage{generic}
\usepackage{cite}
\usepackage{amsmath,amssymb,amsfonts}
\usepackage{algorithmic}
\usepackage{graphicx}
\usepackage{textcomp}

\def\Re{\mathbb{R}}
\def\R{\mathbb{R}}

\def\Theorem#1{Thm.~\ref{#1}}
\def\Corollary#1{Cor.~\ref{#1}}
\def\Sec#1{Sec.~\ref{#1}}

\def\Appendix#1{Appendix~\ref{#1}}
\def\notes#1{\marginpar{\tiny #1}\typeout{Notes!
Notes!
Notes!
}}
\renewcommand{\notes}[1]{\typeout{notes!}}
\def\FRAC#1#2#3{\genfrac{}{}{}{#1}{#2}{#3}}
\def\half{{\mathchoice{\FRAC{1}{1}{2}}%
{\FRAC{2}{1}{2}}%
{\FRAC{3}{1}{2}}%
{\FRAC{4}{1}{2}}}}

\def\Re{\field{R}}

\def\Sec#1{Sec.~\ref{#1}}

\def\clB{{\cal B}}

\def\clL{{\cal L}}
\def\clP{{\cal P}}
\def\clZ{{\cal Z}}

\def\Sec#1{Sec~\ref{#1}}

\def\E{{\sf E}}

\def\clC{{\cal C}}

\def\R{\mathbb{R}}

\def\Sec#1{Sec.~\ref{#1}}

\def\clZ{{\cal Z}}

\newtheorem{theorem}{Theorem}

\newtheorem{definition}{Definition}

\newtheorem{remark}{Remark}
\newtheorem{proposition}{Proposition}
\newtheorem{corollary}{Corollary}

\def\beq{\begin{eqnarray}} 
\def\bc{\begin{center}} 
\def\be{\begin{enumerate}}
\def\bi{\begin{itemize}} 
\def\bs{\begin{small}}
\def\bS{\begin{slide}}
\def\ec{\end{center}} 
\def\ee{\end{enumerate}}
\def\ei{\end{itemize}}
\def\es{\end{small}}
\def\eS{\end{slide}}
\def\eeq{\end{eqnarray}}

\newcommand{\newP}[1]{\medskip\noindent{\bf #1:}}

\newcommand{\ud}{\,\mathrm{d}}

\def\Re{\mathbb{R}}
\def\E{{\sf E}}

\def\clY{{\cal Y}}

\def\Sec#1{Sec.~\ref{#1}}
\def\Thm#1{Thm.~\ref{#1}}
\def\Prop#1{Prop.~\ref{#1}}

\def\clB{{\cal B}}

\def\clL{{\cal L}}
\def\clP{{\cal P}}
\def\clZ{{\cal Z}}

\renewcommand{\Re}{\mathbb{R}}

\def\FRAC#1#2#3{\genfrac{}{}{}{#1}{#2}{#3}}

\def\clA{{\cal A}}
\def\clB{{\cal B}}
\def\clC{{\cal C}}

\def\clF{{\cal F}}

\def\clL{{\cal L}}

\def\clM{{\cal M}}
\def\clN{{\cal N}}
\def\clO{{\cal O}}
\def\clP{{\cal P}}

\def\clU{{\cal U}}

\def\clY{{\cal Y}}
\def\clZ{{\cal Z}}
\def\E{{\sf E}}
\def\bS{\mathbb{S}}

\def\ones{{\sf 1}}
\def\sP{{\sf P}}
\def\tsP{{\tilde{\sf P}}}
\def\tE{{\tilde{\sf E}}}

\def\tp{{\hbox{\rm\tiny T}}}

\def\Nsp{{\sf N}} 
\def\Rsp{{\sf R}} 
\def\hE{\tilde{\sf E}}

\def\dv{\operatorname{diag}}
\def\sp{\operatorname{span}}

\def\bmu{{\bar\mu}}

\def\kl{{\sf D}}

\def\sW{{\sf W}}

\begin{document}
\title{Duality for Nonlinear Filtering I: Observability}
\author{Jin W. Kim, \IEEEmembership{Student member, IEEE}, and
  Prashant G. Mehta, \IEEEmembership{Senior member, IEEE}
\thanks{This work is supported in part by the NSF award 1761622.}
\thanks{Research reported in this paper was carried out by J. W. Kim, as part of his PhD
  dissertation work, while he was a graduate student at the University of Illinois at Urbana-Champaign.  
He is now with the Institute of Mathematics at the University of Potsdam
(e-mail: jin.won.kim@uni-potsdam.de).}
\thanks{P. G. Mehta is with 
the Coordinated Science Laboratory and the Department of Mechanical Science
and Engineering at the University of Illinois at Urbana-Champaign
(e-mail: mehtapg@illinois.edu).}}

\maketitle

\begin{abstract}
This paper is concerned with the development and use of duality theory
for a hidden Markov model (HMM) with white noise
observations. The main contribution of this work is to introduce a
backward stochastic differential equation (BSDE) as a dual control
system.  A key outcome is that stochastic observability
(resp. detectability) of the HMM is expressed in dual terms: as
controllability (resp. stabilizability) of the dual control system.
All aspects of controllability, namely, definition of controllable
space and controllability gramian, along with their properties and
explicit formulae, are discussed.  The proposed duality is 
shown to be an exact extension of the classical duality in 
linear systems theory.  One can then relate and compare the linear and
the nonlinear systems.  A side-by-side summary of this relationship is
given in a tabular form (Table~\ref{tb:comparison}).
\end{abstract}


\section{Introduction}
\label{sec:introduction}

There is a fundamental dual relationship between estimation and
control. The dual relationship is expressed in two inter-related
manners:
\begin{itemize}
	\item Duality between observability and controllability.
	\item Duality between optimal filtering and optimal control. 
\end{itemize}
The second bullet means expressing one type of problem as another type
of problem.  In this two-part paper, the main interest is to
convert a filtering problem into a control problem.

Duality is coeval with the origin of modern systems and control
theory: Early papers of Kalman include a statement of the original
duality principle~\cite[Eq.~(62) and~(72)]{kalman1960general}.  The basic 
papers~\cite{kalman1960,kalman1961} on the Kalman filter
contain an extensive mention of duality 
-- between the optimal filter and a certain linear quadratic optimal control problem.  Notably,
duality explains why the Riccati equation is the fundamental equation
for \textit{both} optimal filtering and optimal control.

Sixty years have elapsed since Kalman's original work. One
would imagine that duality for the nonlinear stochastic systems
(hidden Markov models) is well understood by now. It is a foundational
question at the heart of modern systems and control theory, and its
modern avatars such as reinforcement learning. However, this is not
the case! In his 2008 paper~\cite{todorov2008general}, Todorov writes:
\begin{quote}
	``\it{Kalman's duality has been known for half a century and has attracted a lot of attention. If a straightforward generalization to non-LQG settings was possible it would have been discovered long ago. Indeed we will now show that Kalman's duality,
	although mathematically sound, is an artifact of the LQG setting.}''
\end{quote}
Is this to suggest that there is no previous work to extend duality to nonlinear systems? Au contraire!
For deterministic systems, almost every definition of nonlinear observability, and there have been several notable ones throughout the
decades, appeals to duality.  Likewise, Mortensen and related minimum
energy algorithms, originally invented in 1960s, are standard
approaches to construct an estimator.  For stochastic systems as well,
there have been seminal contributions, notably the work of
Mitter and co-authors~\cite{fleming1982optimal,mitter2003}. 
Having said that, several reasons are noted
in~\cite{todorov2008general} on why the
duality described in these prior works are {\em not} generalizations
of the original Kalman-Bucy duality.  A comprehensive account of the
differences is contained in~\cite[Ch.~3]{JinPhDthesis}.

\subsection{Summary of original contributions}
We consider the stochastic filtering problem
for a hidden Markov model (HMM) with white noise observations (the
mathematical model is introduced in
\Sec{sec:problem-formulation}). For this filtering problem, we
make two types of original contributions:
\begin{itemize}
	\item Dual controllability characterization of stochastic
          observability.  It is the subject of this present paper
          (part I).
	\item Dual (minimum variance) optimal control formulation of the
          stochastic filtering problem. It is the subject of a
          companion paper (part II) also submitted to this journal.
\end{itemize}

The focus of this paper is on the dual controllability characterization of stochastic observability of an HMM with white noise observations.
The dual control system is introduced in
\Sec{sec:dual-control-system}.  Once the dual control system has been
introduced, the ensuing considerations are entirely parallel to linear
systems theory: 
The solution operator of the dual control system is used to define a
linear operator whose range space is the controllable subspace.  The
system is said to be controllable if the range space is dense in a suitable
function space.  The
controllability (resp. stabilizability) of the dual system is shown to be equivalent to
stochastic observability (resp. detectability) of the HMM. 
Several properties of the
controllable subspace are noted along with its explicit
characterization in the finite state-space case.  A formula for the
controllability gramian is also described.   
An upshot of our work is that we can establish parallels between
linear and nonlinear models (Table~\ref{tb:comparison}).    

The part II builds on the results of the part I.  In particular, duality is used to transform the minimum variance objective of the nonlinear filtering problem into a
stochastic optimal control problem.  For the latter problem, the dual
control system, introduced in part I, is shown to arise as the constraint.


\subsection{Relationship to literature}

In linear systems theory, the following systems are said to be dual to each
other:
\begin{align}
\text{(state-output)} \quad\; \dot{x}_t &= A^\tp  x_t,\;x_0=\xi \label{eq:LTI-obs-intro} \\
z_t &= H^\tp x_t,\qquad\qquad 0\le t\le T \nonumber\\
\text{(state-input)} \;\;-\dot{y}_t &= A y_t + H u_t,\; y_T = 0, \;\;0\leq t\leq T \label{eq:LTI-ctrl-intro}
\end{align}
The states, $x_t$ and $y_t$ for the two systems, are vector-valued,
both of
dimension $d$ (the standard dot product in $\Re^d$ is denoted by
$\langle \cdot, \cdot \rangle_{\Re^d}$).  The input $u:=\{u_t\in\Re^m:0\leq t\leq T\}$
and the output $z:=\{z_t \in\Re^m:0\leq t\leq T\}$ are elements of the
function space $L^2([0,T];\Re^m)=:\clU$ equipped with the
inner-product $\langle u, v \rangle_\clU:=\int_0^T u_t^\tp v_t \ud t$. 
For the state-input system~\eqref{eq:LTI-ctrl-intro}, the solution
map $u\mapsto y_0$ is used to define a linear operator $\clL: \clU \to \Re^d$ as follows:
\begin{equation*}\label{eq:LTI-adjoint}
	\clL u := y_0 = \int_0^T e^{At}Hu_t \ud t
\end{equation*}
Its adjoint is given by
\[
(\clL^\dagger \xi)(t) = H^\tp e^{A^\tp t}\xi,\quad 0\le t \le T
\]
and represents the solution map from the initial condition $\xi
\mapsto z$ for the state-output
system~\eqref{eq:LTI-obs-intro}.  Mathematically, the 
dual relationship is expressed as
\[
\langle \xi,\clL u\rangle_{\Re^d} = \langle \clL^\dagger \xi, u \rangle_\clU,\quad \forall\, \xi \in \Re^d,\; u\in \clU
\]
The relationship yields the following
important identity (also referred to as the closed range theorem)
\[
\Rsp(\clL)^\bot = \Nsp(\clL^\dagger)
\]
This identity has several important consequences, e.g., controllability
(resp. stabilizability) 
property of the state-input system~\eqref{eq:LTI-ctrl-intro} is equivalent to the
observability (resp. detectability) property of the state-output system~\eqref{eq:LTI-obs-intro}.

This classical duality between controllability and observability is
useful for \emph{both} analysis and the design of estimation
algorithms.  
For example, most proofs of stability of the Kalman
filter (see e.g.,~\cite[Ch.~9]{xiong2008introduction}) rely -- in
direct or indirect fashion -- on duality theory.  Specifically~(i)
Because of duality, asymptotic stability of the Kalman filter is equivalent to 
asymptotic stability of the (dual) optimal control system; 
(ii)~necessary and
sufficient conditions for the same are stabilizability for the 
control problem, and (because of duality) detectability for the
estimation problem; 
and (iii)~analysis of the optimal control problem (e.g., convergence of
the value function to its stationary limit) yields useful conclusions
on asymptotic filter stability.

This has naturally spurred a large body of work related to:
\begin{itemize}
\item Defining observability as a dual property.
\item Using these definitions to investigate asymptotic 
stability of optimal and sub-optimal estimators.
\end{itemize}
The second bullet has by far been the most important reason for
defining and studying observability and related concepts.  In many
studies, the definition of observability is often a sufficient
condition that guarantees the stability of the estimator under study, e.g.,
~\cite[Defn.~2]{krener2003convergence},~\cite[Defn.~4.1]{rawlings2017model},~\cite[Assumption
2]{copp2017simultaneous},~\cite[identifiability conditions A-1 and
A-2]{baxendale2004asymptotic}.

We next provide a brief survey of observability and its use for
investigating estimator (resp. filter) stability, in the study of
deterministic (resp. stochastic) nonlinear models.  
These have to be separated because the
work on these two types of models has little overlap. For
deterministic models, an estimator is defined as a dynamical system
whose dimension is the same as the dimension of the state and which 
operates on inputs and outputs.  An important class, which also incorporates certain optimality properties, is the
minimum energy estimator
(MEE)~\cite[Ch.~4]{rawlings2017model}. (Relationship to MEE
is described in part II.)   
The two bullets above serve as guiding
principles around which the discussion is organized.

\medskip

\noindent \textbf{Deterministic models.} 
In the classical
paper~\cite{hermann1977nonlinear}, Hermann and Krener write ``{\it
  duality between ``controllability'' and ``observability'' [...] is,
  mathematically, just the duality between vector fields and
  differential forms}''.  Similarly, the output-to-state stability
(OSS) definitions in~\cite{sontag1997output} are
motivated by Wang and Sontag as follows:
``{\it Given the central role often played in control theory by
  the duality between input/state and state/output behavior, one may
  reasonably ask what concept obtains if outputs are used instead of
  inputs in the [input-to-state stability (ISS)] definition''}.  The
OSS definition is important for the following two
reasons: (i) For the linear state-output
system~\eqref{eq:LTI-obs-intro}, OSS is equivalent to
detectability~\cite[Excercise 7.3.12]{sontag2013mathematical}; and
(ii) OSS admits a dissipative
characterization in terms of an OSS-storage
function~\cite[Thm.~3]{sontag1997output}.  Such characterizations are
useful in the study of stability and robustness (several variations of observability
definition and their relationship are discussed
in~\cite{sontag2008input,hespanha2002nonlinear}). 
Combining ISS and OSS yields the notion of
IOSS which is shown to be equivalent to estimating the norm of
the hidden state~\cite[Sec.~8.5]{sontag2008input}.  Related to
detectability, an important notion is the incremental IOSS (i-IOSS)
which is standard for asymptotic stability analysis of MEE ~\cite[Thm.~4.10]{rawlings2017model}.  Dissipative characterizations of
incremental notions are also important, e.g., an i-IOSS Lyapunov function is given
in~\cite{allan2021nonlinear}.  

\medskip

\noindent \textbf{Stochastic models (HMMs).}  As in the study of
deterministic models, a major impetus to define
observability/detectability comes from the question of nonlinear
filter stability (now in the sense of asymptotic forgetting of the initial
prior).  Formally, there are two main cases~\cite{mcdonald2022robustness}: 
\begin{itemize}
\item The case where the Markov process forgets
the prior and therefore the filter ``inherits'' the same property; 
\item The case where the observation
provides sufficient information about the hidden state, allowing the
filter to correct its erroneous initialization.
\end{itemize}
These two cases are referred to as the ergodic and non-ergodic signal
cases, respectively.  While the two cases are intuitively reasonable,
they spurred much work during 1990-2010 with a complete
resolution appearing only at the end of this time-period.  For the ergodic
case, sufficient conditions for filter stability that rely only on the
signal model appear in~\cite[Thm.~5]{atar1997lyapunov},~\cite[Assumption
4.3.24]{Moulines2006inference},~\cite[Thm.~4.3]{baxendale2004asymptotic},~\cite[Cor.~2.3.2]{van2006filtering},~\cite[Thm.~4.2]{baxendale2004asymptotic}. For
the non-ergodic signal case, sufficient conditions relying
also on the observation model are given in~\cite[Thm~.7]{atar1997lyapunov},~\cite[A-1 and
A-2]{baxendale2004asymptotic},~\cite[Rem.~23.1]{stannat2006}. 
All of these conditions are for the HMM with white noise observations,
a model which occupies a central place in the nonlinear filtering theory.  For a more 
general class of HMMs, the fundamental definition for stochastic
observability and detectability is
due to van Handel~\cite{van2009observability,van2009uniform} (see~\Sec{sec:van-Handel-stoch-obsvbl}).  There
are two notable features: (i) the definition made
rigorous the intuition described in the two
cases~\cite[Sec. II-B and Sec. V]{van2010nonlinear}; and (ii) the
definition led to meaningful conditions that were shown to be 
necessary and sufficient for filter stability~\cite[Thm.~III.3 and
Thm.~V.2]{van2010nonlinear}. The 
stochastic observability definition is entirely probabilistic and its information-theoretic
extension was given in~\cite{liu2011thesis,liu2011stochastic}.  For linear
stochastic systems, information-theoretic metrics such as relative
entropy had earlier been used to define
observability~\cite{ugrinovskii2003observability} (see also~\cite{petersen2002notions} where extensions for
uncertain linear systems is described).   
There are also a number of works where observability is defined as
a finite memory property (often referred to as uniform observability or reconstructability) whereby only the most recent window of
observations is necessary for estimation and/or
control ~\cite[Assumption 2]{copp2017simultaneous},~\cite[Defn.~2.4]{rao2003constrained}~\cite[Assumption A2]{michalska1995moving}.  Similar
considerations also inform~\cite{mcdonald2019cdc,mcdonald2022robustness} where $N$ step
observability is defined for an HMM.  The definition has many
attractive features, e.g., it is relatively easier to compute (as compared
to stochastic observability) and is
useful for filter and control design.
In~\cite[Fig.~1]{mcdonald2022robustness}, the definition is also related
to several criteria for filter stability.

In contrast to the deterministic
models, duality is conspicuous by its absence  both in defining and
in using stochastic observability.  Dissipative characterizations -- that
are so familiar in the study of deterministic models -- are also missing.

This paper is drawn from the PhD thesis~\cite{JinPhDthesis} of the first
author.  A prior conference version of this paper appeared
in~\cite{kim2019observability}.  While the focus of conference paper
was on the finite state-space case, the present journal version includes the
results for the general case.  
The following additions are noted: A novel extension to stabilizability of 
the dual control system is described (\Sec{sec:stab-detec}) and 
shown to be the dual to detectability of the HMM (\Corollary{cor:detect-stab}). Remarks~\ref{rm:terminal-condition},~\ref{rm:duality-pairing}, and~\ref{rm:range-dense} are included
to clarify the choice of function spaces.  The dual control system described
in our work is compared and contrasted to the backward Zakai
equation which is how duality is understood in the theory of nonlinear filtering
(Rem.~\ref{rm:backward-Zakai}).  Finally, the results are related
to both stochastic observability (\Theorem{thm:observability-Zakai-relation}) and to the linear Gaussian problem (\Sec{rm:duality-reduces-LG}).






\subsection{Paper outline}

The outline of the remainder of this paper is as follows: The problem
formulation and background appears in
\Sec{sec:problem-formulation}. 
The dual
control system is introduced in \Sec{sec:dual-control-system} together with the definition of the controllability and related concepts.  The explicit formulae for the
finite state space case appear in
\Sec{sec:finite-observability}.  The paper closes with some
conclusions and directions for future research in
\Sec{sec:conc}.


\section{Background and Problem Formulation}
\label{sec:problem-formulation}

\subsection{Notation}

For a locally compact Polish space $S$, the following notation is adopted:
\begin{itemize}
	\item $\clB(S)$ is the Borel $\sigma$-algebra on $S$.
	\item $\clM(S)$ is the space of regular, bounded and finitely additive 
	signed measures (rba measures) on $\clB(S)$. The natural norm
        is the total variation norm denoted $\|\cdot\|_{\text{\tiny
            TV}}$. 
	\item $\clP(S)$ is the subset of $\clM(S)$ comprising of probability 
	measures.
	\item $C_b(S)$ is the space of continuous and bounded
          real-valued functions on $S$. The natural norm
        is the sup-norm $\|\cdot\|_{\infty}$. 
	\item For measure space $(S;\clB(S);\lambda)$, $L^2(\lambda)= L^2(S;\clB(S);\lambda)$ is the Hilbert space of real-valued functions equipped with the inner product
$	
	\langle f, g\rangle_{L^2(\lambda)} = \int_S f(x)g(x)\ud \lambda(x)
$. 
\end{itemize}
For functions $f:S\to \Re$ and $g:S\to \Re$, the notation $fg$ is used to denote element-wise 
product of $f$ and $g$, namely,
$
(fg)(x) := f(x)g(x)
$ for $x\in S$. 
In particular, $f^2 = ff$. The constant function is denoted by $\ones$ ($\ones(x) = 1$ for all $x\in S$). 
For $\mu\in \clM(S)$ and $f\in C_b(S)$, $
\mu(f) := \int_S f(x) \ud \mu(x) $ and for $\mu,\nu\in \clM(S)$ such that $\mu$ is absolutely continuous with 
respect to $\nu$ (denoted $\mu\ll\nu$), the Radon-Nikodym (RN) derivative is 
denoted by $\dfrac{\ud \mu}{\ud \nu}$.   For a subset $B\subset
C_b(S)$, the annihilator of $B$, denoted by $B^\bot$, is defined by $
B^\bot := \{\mu\in\clM(S): \mu(f) = 0\; \forall f\in B\}$.   
For a sub $\sigma$-algebra
${\cal G} \subset \clB(S)$, the restriction of the measure $\mu$ to ${\cal
  G}$ is denoted by $\mu|_{\cal G}$.  It is obtained from the
defining relation $\mu|_{\cal G}(B)=\mu(B)$ for
$B\in {\cal G}$.



\subsection{Hidden Markov Model}

We consider continuous-time
stochastic processes on a finite time-horizon $[0,T]$ with a fixed
$T<\infty$. Fix the probability space $(\Omega, \clF_T, \sP)$ along
with the filtration $\{\clF_t:0\le t \le T\}$ with respect to which all the stochastic
processes are adapted. Of special interest is the pair $(X,Z)$ defined as follows:
\begin{itemize}
	\item The \emph{state process} $X = \{X_t:\Omega\to \bS:0\le t \le
          T\}$ is a Feller-Markov 
	process on the state-space $\bS$. Its initial measure (prior) is denoted by 
	$\mu \in \clP(\bS)$ and $X_0\sim \mu$. The infinitesimal generator is denoted by $\clA$.
	\item  The \emph{observation process} $Z = \{Z_t:0\le t \le T\}$ satisfies the stochastic differential equation (SDE):
	\begin{equation*}\label{eq:obs-model}
		Z_t = \int_0^t h(X_s) \ud s + W_t,\quad 0\le t \le T
	\end{equation*}
	where $h:\bS\to \Re^m$ is the observation function and $W = \{W_t:0\le t \le T\}$ is an $m$-dimensional Brownian motion (B.M.). We write $W$ is $\sP$-B.M. 
	It is assumed that $W$ is independent of $X$.
\end{itemize}
The above is referred to as the \emph{white noise observation model} of 
nonlinear filtering. In the remainder of this paper, the model is
referred to as the HMM $(\clA,h)$.  In the case where $\bS$ is not finite, additional assumptions are 
typically necessary to ensure that the model is well-posed.  In
applications, the 
most important examples are as follows:
\begin{itemize}
\item $\bS$ is finite with cardinality
$|\bS|=d$.
\item $\bS=\Re^d$ and $X$ is an
It\^{o} diffusion.  
\end{itemize}
A historically noteworthy example of an It\^{o} diffusion is the linear
Gaussian model. 

\subsection{Nonlinear filtering background}

The canonical filtration $\clF_t = \sigma\big(\{(X_s,W_s):0\le s \le t\}\big)$. The filtration generated by the observation is denoted by $\clZ :=\{\clZ_t:0\le t\le T\}$  where $\clZ_t = \sigma\big(\{Z_s:0\le s\le t\}\big)$. 
The \emph{nonlinear (or stochastic) filtering problem} is to compute the conditional expectation for a given test function $f\in C_b(\bS)$:
\[
\pi_t(f) := \E\big(f(X_t)\mid \clZ_t\big),\quad 0\le t \le T
\]
The measure-valued process $\pi = \{\pi_t:0\le t \le T\}$ is referred to as the nonlinear filter.

A standard approach 
is based upon the Girsanov change of measure.
Suppose the model satisfies the Novikov's condition: $
\E\left(\exp\big(\half \int_0^T |h(X_t)|^2\ud t\big)\right) < \infty
$. 
Define a new measure $\tsP$
on $(\Omega,\clF_T)$ as follows:
\[
\frac{\ud \tsP}{\ud \sP} = \exp\Big(-\int_0^T 
h^\tp(X_t) \ud W_t - \half \int_0^T |h(X_t)|^2\ud t\Big) =: D_T^{-1}
\]
Then it can be shown that the probability law for $X$ is unchanged but
$Z$ is a $\tsP$-B.M.~that is independent of $X$~\cite[Lem.~1.1.5]{van2006filtering}. The
expectation with respect to $\tsP$ is denoted by $\tE(\cdot)$.  
The \emph{un-normalized filter} is a 
measure-valued process $\sigma = \{\sigma_t:0\le t \le
T\}$ defined by
\[
\sigma_t(f) := \tE\big(D_tf(X_t)|\clZ_t\big),\quad f\in 
C_b(\bS)
\]
Because $Z$ is a $\tsP$-B.M., the equation is unnormalized filter is
easily obtained and is in fact the celebrated Zakai equation of
nonlinear filtering \cite[Thm.~5.5]{xiong2008introduction}:
	\begin{equation}\label{eq:Zakai}
		\sigma_t(f) = \mu(f) + \int_0^t \sigma_s(hf)^\tp \ud Z_t + 
		\int_0^t \sigma_s(\clA f)\ud s,\quad 0\leq t\leq T
	\end{equation}
Upon normalization, the nonlinear filter
\begin{equation}\label{eq:normalize-Zakai}
	\pi_t(f) = \frac{\sigma_t(f)}{\sigma_t(\ones)},\quad 0\leq t\leq T
\end{equation}
This ratio is referred to as the Kallianpur-Striebel formula~\cite[Thm.~5.3]{xiong2008introduction}. Using~\eqref{eq:Zakai}, the equation for the nonlinear
filter is readily obtained by a simple application of the It\^o formula to the ratio~\cite[Thm.~5.7]{xiong2008introduction}.

\subsection{Stochastic
  observability}\label{sec:van-Handel-stoch-obsvbl}

In problems concerned with observability of the model $(\clA,h)$ or
filter stability, there are reasons to consider more than one prior. We reserve the notation $\mu$ to denote the true but possibly
unknown prior and the notation $\nu$ to denote the prior that is used
to compute the filter. If $\mu$ is exactly known then $\mu = \nu$. In
all other cases, it is assumed that $\mu\ll\nu$.

To stress the dependence on the initial measure $\mu$, we use the superscript 
notation $\sP^\mu$ to denote the probability measure $\sP$ when $X_0\sim \mu$. 
The expectation operator is denoted by $\E^\mu(\cdot)$ and the nonlinear filter 
$\pi_t^\mu(f) = \E^\mu\big(f(X_t)| \clZ_t\big)$. On the common measurable 
space $(\Omega, \clF_T)$, $\sP^\nu$ is used to denote another probability measure 
such that the transition law of $(X,Z)$ are identical but $X_0\sim \nu$. (For an explicit construction of $\sP^\mu$ and $\sP^\nu$, see~\cite[Sec.~2.2]{clark1999relative}.)  The 
associated expectation operator is denoted by $\E^\nu(\cdot)$ and 
$\pi_t^\nu(f) = \E^\nu\big(f(X_t)|\clZ_t\big)$.  The respective un-normalized
filters are denoted by $\sigma_t^\mu$ and $\sigma_t^\nu$.  These are
solution of the Zakai equation~\eqref{eq:Zakai} with initialization
$\sigma_0^\mu=\mu$ and  $\sigma_0^\nu=\nu$, respectively.

The following definition of stochastic observability is introduced in~\cite{van2009observability}.
Although it is identical for a general class of HMMs, we state it for the model $(\clA,h)$:

\medskip

\begin{definition}[Defn.~2 in~\cite{van2009observability}]\label{def:observability}
	The model $(\clA,h)$ is \emph{observable} if 
\begin{equation*}
	\quad \sP^\mu|_{\clZ_T} = 
	\sP^\nu|_{\clZ_T} \; \Longrightarrow\; \mu = \nu,\quad \forall\,\mu,\nu\in\clP(\bS)
\end{equation*}
\end{definition} 

\medskip

\begin{remark}
The definition is contrasted with the definition of observability for
deterministic nonlinear models~\cite[Defn.~6.1.4]{sontag2013mathematical}.  For deterministic models, observability is
defined as a property of the map from initial condition to output
trajectory.  In contrast, stochastic observability is a property of
the map from the initial prior to the probability law of the output
process. 
\end{remark}

\medskip

Consider an equivalence relation on $\clP(\bS)$ as follows: 
\[
\mu \simeq \nu\quad \text{if} \quad \sP^\mu|_{\clZ_T} = 
\sP^\nu|_{\clZ_T}
\]
The following definition naturally arises from this notation:

\medskip

\begin{definition}[Defn.~3 
in~\cite{van2009observability}]\label{def:un-observable-measures-observable-functions}
	The space of \emph{observable functions} 
$
	\clO = \{f\in C_b(\bS): \mu(f) = \nu(f) \; \forall\, \mu\simeq
        \nu\}$. 
	The space of \emph{unobservable measures} 
$
	\clN = \{c(\mu - \nu) \in \clM(\bS): c \in \Re,\; \mu,\nu 
	\in \clP(\bS)\; \text{s.t.}\; \mu\simeq\nu\}
$.
\end{definition} 


\medskip

\begin{remark}
An HMM is observable if and only if $\clN =
\{0\}$, and because $\clO^\bot = \clN$, an HMM is observable
if and only if $\clO$ is dense in
$C_b(\bS)$~\cite[p.~42]{van2009observability}. Because $\ones\in\clO$,
the space of observable functions $\clO$ is always non-trivial.  This
also means $\clN\subseteq \clM_0(\bS):=\{\mu\in\clM(\bS): \mu(\ones) = 0\}$.
\end{remark}

\section{Dual Control System}
\label{sec:dual-control-system}

\subsection{Function spaces}

It is noted that $Z$ is a $\tsP$-B.M..  
For a $\clZ_T$-measurable
random vector, the following definition of Hilbert space is
standard: $
L^2_{\clZ_T}(\Omega;\Re^m) := L^2(\Omega;\clZ_T;\ud
\tsP)$~\cite[Ch.~5.1.1]{le2016brownian}.  
For a $\clZ$-adapted vector-valued stochastic process, the Hilbert space is
$L^2_{\clZ}\big([0,T];\Re^m\big):=
L^2\big(\Omega\times[0,T];\clZ\otimes \clB([0,T]);\ud \tsP\ud
t\big)$ 
where $\clB([0,T])$ is the Borel sigma-algebra on $[0,T]$,
$\clZ\otimes \clB([0,T])$ is the product sigma-algebra and $\ud
\tsP\ud t$ denotes the product measure on it. The inner product for these spaces are
\[
\langle F,G\rangle_{L^2_{\clZ_T}} = \tE\big(F^\tp G\big),\quad \langle U,V\rangle_{L^2_\clZ} = \tE\Big(\int_0^T U_t^\tp V_t \ud t\Big)
\]
These Hilbert spaces suffice if the state-space $\bS$ is finite.  In general
settings, let $\clY$ denote a suitable Banach space of real-valued
functions on $\bS$, equipped with the norm $\|\cdot\|_\clY$. 
Then
\begin{itemize}
\item For a random function, the
Banach space $L^2_{\clZ_T}(\Omega;\clY) := \big\{F:\Omega\to
\clY: F\text{ is }\clZ_T\text{-measurable},\;
\tE\big(\|F\|_\clY^2\big) < \infty\big\}$.  
\item 
For a function-valued
stochastic process, the Banach space is
$L^2_{\clZ}([0,T];\clY) := \Big\{Y:\Omega\times
[0,T]\to \clY\;:\; Y\text{ is }\clZ\text{-adapted},\;
\tE\Big(\int_0^T\|Y_t\|_\clY^2 \ud t\Big) < \infty\Big\}$. 
\end{itemize}
In this paper, examples of $\clY$ are: (i) $C_b(\bS)$ equipped with
sup norm denoted by $\|\cdot\|_\infty$, and (ii) $L^2(\lambda)$ where
$\lambda$ is a positive reference measure on $\bS$. 


\subsection{Stochastic observability and Zakai equation}

For a white noise observation model, a quantitative analysis of
stochastic observability is possible based on the following formula
for relative entropy~\cite[Thm.~3.1]{clark1999relative}:
	\begin{equation}\label{eq:clark-result}
	\kl\big(\sP^\mu|_{\clZ_T} \mid \sP^\nu|_{\clZ_T}\big) = \half\E^\mu\Big(\int_0^T |\pi_t^\mu(h)-\pi_t^\nu(h)|^2\ud t\Big)
	\end{equation}
where $\kl(\cdot \mid \cdot)$ is the Kullback-Leibler (KL) divergence.  The
formula is used to obtain the following result: 

\medskip 

\begin{theorem}\label{thm:observability-Zakai-relation}
	T.F.A.E.:
	\begin{enumerate}
		\item The model $(\clA, h)$ is observable.
		\item For $\mu,\nu\in\clP(\bS)$,
		\[
		\pi_t^\mu(h) = \pi_t^\nu(h),\quad t\text{-a.e.},\;\sP^\mu|_{\clZ_T}\text{-a.s.} \quad \Longrightarrow \quad \mu = \nu
		\]
		\item For $\mu,\nu\in\clP(\bS)$,
		\[
		\sigma_t^\mu(h) = \sigma_t^\nu(h),\quad t\text{-a.e.},\; \sP^\mu|_{\clZ_T}\text{-a.s.} \quad \Longrightarrow \quad \mu = \nu
		\]
	\end{enumerate}
\end{theorem}

\medskip

\begin{proof}
See \Appendix{pf-obs-Zakai}.
\end{proof}

\medskip

The value of Thm.~\ref{thm:observability-Zakai-relation} is that the un-normalized filter is the solution to the Zakai equation which is linear.
A linear operator $\clL^\dagger: \clM(\bS)\to L^2_\clZ\big([0,T];\Re^m\big)\times \Re$ is defined as follows:
\begin{equation}\label{eq:observability-operator}
	\clL^\dagger \mu = \big(\{\sigma_t^\mu(h):0\le t\le T\}, \mu(\ones)\big)
\end{equation}
The notation is suggestive: In this paper, we will define a linear operator $\clL$ such that the operator defined by~\eqref{eq:observability-operator} is its adjoint.

\medskip

\begin{corollary}\label{cor:linear-operator-observability}
	The model $(\clA,h)$ is observable if and only if $
	\Nsp(\clL^\dagger) = \{0\}$. 
\end{corollary}

\medskip

\begin{remark}
	$\tilde{\mu}\in\Nsp(\clL^\dagger)$ is the space of unobservable measure (see Defn.~\ref{def:un-observable-measures-observable-functions}). Suppose
        $\Nsp(\clL^\dagger)$ is non-trivial. Then for
        $\mu\in\clP(\bS)$, choose $\epsilon \neq 0$ such that $\nu =
        \mu+\epsilon \tilde{\mu} \in \clP(\bS)$. Then owing to the
        linearity of~\eqref{eq:Zakai}, $
	\sigma_t^\mu(h) = \sigma_t^\nu(h)$ for $0\leq t\leq T$.  
	From \Theorem{thm:observability-Zakai-relation} then
        $\sP^\mu|_{\clZ_T} = \sP^\nu|_{\clZ_T}$.  
\end{remark}

\medskip

\begin{remark}[Prior work on adjoint of the Zakai equation]  \label{rm:backward-Zakai}
	Because the Zakai equation is linear, its adjoint has previously been
	considered in literature. There are two types of equivalent constructions:
	\begin{itemize}
		\item The most direct route relies on the pathwise or the robust
		representation of the nonlinear filter.  In this approach, by using
		a log transformation, the stochastic partial differential
		equation (PDE) is transformed into a linear deterministic PDE with
		random coefficients~\cite[Ch.~VI.11]{rogers2000diffusions}.  Because
		the transformed PDE is deterministic, the formula for adjoint is obtained by
		standard means~\cite[Eq.~4.17-4.18]{benevs1983relation}.    
		\item 
		The second type of adjoint is the \emph{backward Zakai equation}
		\begin{equation}\label{eq:backward-Zakai-eqn}
			\begin{aligned}
				-\ud \eta_t(x) & = (\clA \eta_t)(x) \ud t
				+ (h\eta_t)(x)
				\overleftarrow{\ud Z_t}\\ \eta_T (x)&=
				f(x), \quad x\in \bS
			\end{aligned}
		\end{equation}
		where $\overleftarrow{\ud Z_t}$ denotes the backward It\^o integral:
		its construction is based on choosing the right-endpoints in the
		partial sum approximation of the stochastic
		integral~\cite[Rem.~3.3]{pardoux1981non}.  The backward and forward Zakai equation are said to be dual
		because of the following formula~\cite[Thm.~4.7.5]{bensoussan1992stochastic}:
		\[
		\sigma_T(f) = \mu(\eta_0)
		\]
		The formula is used to prove the uniqueness of the solution to the
		(forward)
		Zakai equation~\cite[Sec.~6.5]{xiong2008introduction}.
	\end{itemize}
	
	The two constructions are equivalent because, using the log
	transformation, the backward Zakai equation is transformed to the
	pathwise adjoint~\cite[Sec.~2.3]{kim2020smoothing}.  In nonlinear
	filtering, the forward and backward Zakai equation are both
	classical~\cite{pardoux1979backward}. The two equations 
	together yields the solution of the smoothing
	problem~\cite[Thm.~3.8]{pardoux1981non}.

	
	Despite the well known duality between forward and backward Zakai equations, it is distinct
	from the controllability--observability duality (in linear systems
	theory) because of the following aspects:
	\begin{itemize}
		\item The dual equation~\eqref{eq:backward-Zakai-eqn} does not have a control input term.
		\item The stochastic process $\eta=\{\eta_t:0\leq t\leq T\}$ is not adapted to any forward-in-time filtration. In particular, $\eta_0$ is a $\clZ_T$-measurable random variable.
	\end{itemize}
	The dual control system described in the following section is original and distinct from these prior adjoints.
	
\end{remark}

\subsection{Dual control system}

The goal is to define a linear operator $\clL$ whose adjoint
$\clL^\dagger$ is given by~\eqref{eq:observability-operator}.  Because
$\clL^\dagger: \clM(\bS)\to L^2_\clZ\big([0,T];\Re^m\big)\times \Re$, 
the domain of $\clL$ is $L_\clZ^2\big([0,T];\Re^m\big)\times \Re$.  We
set $\clU:=L_\clZ^2\big([0,T];\Re^m\big)$ and refer to it as the
space of {\em admissible control inputs}.  Next, because of duality
pairing between $C_b(\bS)$ and $\clM(\bS)$, the co-domain of $\clL$ is
$C_b(\bS)$.  We set $\clY=C_b(\bS)$.

The main result (\Thm{thm:observability-definition} below) is to show that the operator $\clL: \clU\times \Re \to \clY$ is defined by the solution operator of the linear backward stochastic differential equation (BSDE):
%
\begin{subequations}\label{eq:dual-bsde}
	\begin{align}
		-\ud Y_t(x) &= \big((\clA Y_t)(x) + h^\tp(x)(U_t+V_t(x))\big)\ud t - V_t^\tp (x) \ud Z_t\label{eq:dual-bsde-a}\\
		Y_T(x) &= c,\quad\forall\,x\in\bS\label{eq:dual-bsde-b}
	\end{align}
\end{subequations}
where the control input $U:=\{U_t: 0\leq t\leq T\}\in \clU$ and $c\in \Re$ is a deterministic constant.  The
solution of the BSDE $(Y,V) :=\{(Y_t,V_t) : 0\leq t\leq T\}\in
L^2_{\clZ}\big([0,T];\clY\times \clY^m\big)$ is
(forward) adapted to the filtration $\clZ$.
Existence, uniqueness, and regularity theory for linear BSDEs is
standard and 
throughout the paper, we assume that the solution of BSDE $(Y,V)$ is
uniquely determined in
$L^2_{\clZ}\big([0,T];\clY\times \clY^m\big)$ for
each given $Y_T\in L^2_{\clZ_T}(\Omega;\clY)$ and $U\in
L_\clZ^2\big([0,T];\Re^m\big)$.  The well-posedness results for finite
state-space can be found in~\cite[Ch.~7]{yong1999stochastic} and for
the Euclidean state space 
in~\cite{ma1999linear} (see also Rem.~\ref{rm:duality-pairing} below). 


The linear operator $\clL:\clU\times \Re \to \clY$ is defined as follows:
\begin{equation}\label{eq:ctrl-operator}
	\clL(U,c) = Y_0
\end{equation}
where $Y_0\in \clY$ is the solution at time 0 to the
BSDE~\eqref{eq:dual-bsde}. Note that $Y_0$ is a deterministic
function.  

Controllability is now defined in the same way as linear systems theory. Note however that the target set (at time $T$) is the space of constant functions (see also Rem.~\ref{rm:terminal-condition}).

\begin{definition}
	For the BSDE~\eqref{eq:dual-bsde}, the \emph{controllable subspace}
	\begin{equation}\label{eq:ctrl-subspace}
		\clC_T := \Rsp(\clL) 
	\end{equation}
Explicitly $\clC_T = \big\{y_0 \in \clY: \exists\, c\in \Re\text{
  and } U \in \clU\; \text{such that}\; Y_0 = y_0\;\text{and}\; Y_T = c\ones\big\}$.
	The BSDE~\eqref{eq:dual-bsde} is said to be \emph{controllable} if $\clC$ is dense in $\clY$.
\end{definition}
\medskip

The duality between observability of the model $(\clA, h)$ and the controllability of the BSDE~\eqref{eq:dual-bsde} is described in the following theorem:

\medskip

\begin{theorem}\label{thm:observability-definition}
$\clL^\dagger$ is the adjoint operator of $\clL$. Consequently, the HMM $(\clA, h)$ is observable if and only if the BSDE~\eqref{eq:dual-bsde} is controllable. 
\end{theorem}

\medskip

\begin{proof}
See \Appendix{ssec:pf-observability-definition}.
\end{proof}

\medskip

We refer to the BSDE~\eqref{eq:dual-bsde} as the \emph{dual control
  system} for the HMM $(\clA, h)$.  
We make some remarks on function and measure spaces.


\medskip

\begin{remark}\label{rm:terminal-condition}
	In the definition of $\clL$, the domain space is
        $L^2_{\clZ}\big([0,T];\Re^m\big)\times \Re$. For the purposes
        of this 
        study, it also suffices to consider a restriction of
        $\clL^\dagger$ on the subspace $\clM_0(\bS)$ (for example, the
        null-space of $\clL^\dagger$ is a subspace of $\clM_0(\bS)$).
        An advantage of considering such a restriction is that the co-domain space for $\clL^\dagger$, and therefore the domain of $\clL$, now is $L^2_{\clZ}\big([0,T];\Re^m\big)$. The dual space of $\clM_0(\bS)$ is the quotient space $C_b(\bS)/\{c\ones:c\in\Re\}$ and therefore
	$\clL:L^2_{\clZ}\big([0,T];\Re^m\big)\to
        C_b(\bS)/\{c\ones:c\in\Re\}$. Although such a change will make
        duality between controllability and observability somewhat
        terser, we prefer to keep the measure space as $\clM(\bS)$ and
        the function space as $C_b(\bS)$. 
	This has an advantage of not having to deal with the solutions
        of the BSDE on the quotient space.
\end{remark}

\medskip

\begin{remark}\label{rm:duality-pairing}
	The choice of function space $\clY=C_b(\bS)$ is guided by duality
        pairing between $C_b(\bS)$ and measure space
        $\clM(\bS)$~\cite[Thm.~IV.6.2]{dunford1958linear}. An important reason is to relate with~\cite{van2009observability} who defines observable
        functions as a subspace of $C_b(\bS)$. However, this may 
        place restriction on the model (e.g., $\bS$ is finite or
        compact) for the linear operator $\clL:\clU\times \Re\to \clY$ to be bounded.  
Alternatively, one may consider linear operator on a suitable Hilbert
space.  An important case is when a $\bS$ admits a positive
reference measure $\lambda$.  
In this case, provided these are well-defined, one may consider
\begin{align*}
	\clL\;:\;\clU \times \Re
        \to L^2(\lambda),
\qquad
\clL^\dagger\;:\; L^2(\lambda) \to \clU\times \Re
\end{align*}
where the domain for the adjoint $\clL^\dagger$ is the space of
absolutely continuous measures $\nu\in{\clM}(\bS)$ whose density $\dfrac{\ud \nu}{\ud
  \lambda} \in L^2(\lambda)$.  Examples are (i)
finite state-space where $\lambda$ is the counting measure; and (ii)
the Euclidean state-space where $\lambda$ is the Lebesgue
measure.  In fact, the well-posedness results for the BSDE in
Euclidean settings are for $\clY=L^2(\lambda)$~\cite{ma1999linear}. 
For linear Gaussian problems, one may take $\lambda$ to be the
Gaussian prior.  
\end{remark}

\medskip

\begin{remark}\label{rm:range-dense}
In finite state-space settings, $C_b(\bS)$ and $\clM(\bS)$ are
isomorphic to $\Re^d$.  In Euclidean settings, the problem is
infinite-dimensional and 
$\Rsp(\clL)$ can never be the entire function space
$C_b(\bS)$~\cite[Ch.~4]{curtain2012introduction}.  Therefore, the
controllability is (necessarily) defined as the property that the range space is
dense in $C_b(\bS)$.
\end{remark}

\medskip

In the remainder of this paper, we let $\clY$ to be a suitable
function space and $\clY^\dagger$ is the dual space.  A reader may
replace $\clY=C_b(\bS)$ and $\clY^\dagger=\clM(\bS)$ and such a choice
is always possible if $\bS$ is finite or compact.

\subsection{Controllable subspace and controllability
  gramian} \label{ssec:gramian}


The following proposition describes an important property of the
controllable subspace which is useful for computations. 
\begin{proposition}\label{thm:controllable-subspace}
The controllable subspace $\clC_T$  is the smallest such subspace $\clC\subset \clY$ that satisfies the
	following two properties:
	\begin{enumerate}
		\item The constant function $\ones\in \clC$; and
		\item If $g\in\clC$ then $\clA g \in \clC$ and $g h
		\in\clC$. 
	\end{enumerate}
\end{proposition}

\medskip

\begin{proof}
See \Appendix{ss:pf-prop43}.
\end{proof}

\medskip

The \emph{controllability gramian} $\sW :\clY^\dagger  \to \clY$ is a deterministic linear operator defined as follows:
\[
\sW := \clL \clL^\dagger
\]	
Explicitly, for a measure $\mu \in \clY^\dagger$, $
\sW \mu = Y_0
$
where $Y_0$ is obtained for solving the BSDE
\begin{align*}
-\ud Y_t(x) &= \big(\clA Y_t(x) + h^\tp(x)(\sigma_t^\mu(h)+V_t(x))\big)\ud
t - V_t^\tp(x)\ud Z_t\\
Y_T(x) &= \mu(\ones),\quad x\in \bS
\end{align*}
The gramian yields an explicit control input as follows:

\medskip

\begin{proposition}\label{prop:gramian-observability}
	Suppose $f\in \Rsp(\sW)$, i.e., there exists $\mu\in \clY^\dagger$ such that $f = \sW \mu$. Then the control
	\[
	U_t = \sigma_t^\mu(h),\quad 0\le t \le T
	\]
	transfers the system~\eqref{eq:dual-bsde} from $Y_T = \mu(\ones)\ones$ to $Y_0 = f$.
	Suppose $\tilde{U}$ is another control which also transfers $Y_T = c\ones$ to $Y_0 = f$ for some $c\in \Re$. Then
\begin{equation}\label{eq:min-norm}
	\tE\Big(\int_0^T|\tilde{U}_t|^2 \ud t \Big) + c^2 \ge \tE\Big(\int_0^T |U_t|^2\ud t \Big) + \big(\mu(\ones)\big)^2
	\end{equation}
\end{proposition}


\medskip

\begin{proof}
See \Appendix{ss:pf-prop44}.
\end{proof}


\begin{table*}[t]
	\caption{Linear operators and their adjoints for deterministic (top) and 
		stochastic (bottom) systems}\label{tab:tab1}
	\hrule
	\hrule
	\begin{flalign*}
		&\text{(state-output system Eq.}~\eqref{eq:LTI-obs-intro})&
		\quad \text{(initial condition at $t=0$)} \quad
		&\xi\in \Re^d \quad \stackrel{{\cal
				L}^\dagger}{\longrightarrow} \quad \{z_t:0\leq
		t\leq T\} \in L^2\big([0,T];\Re^m\big) \quad &\text{(output)} && \\
		&\text{(state-input system Eq.}~\eqref{eq:LTI-ctrl-intro})&
		\quad \text{(initial condition at $t=0$)} \quad
		&y_0\in \Re^d \quad \stackrel{{\cal L}}{\longleftarrow} \quad \{u_t :0\leq
		t\leq T\} \in L^2\big([0,T];\Re^m\big)\quad  &\text{(control input)} &&
	\end{flalign*}
\hrule
	\begin{flalign*}
		&\text{(Zakai Eq.}~\eqref{eq:Zakai})&
		\quad \text{(measure at $t=0$)} \quad &\mu\in {\cal
			M}(\mathbb{S}) \quad \stackrel{{\cal
				L}^\dagger}{\longrightarrow} \quad \{\sigma_t^\mu(h):0\leq
		t\leq T\} \in L^2_\clZ\big([0,T];\Re^m\big) \quad &\text{(un-norm. filter)} && \\
		&\text{(BSDE Eq.}~\eqref{eq:dual-bsde})&
		\quad \text{(function at $t=0$)} \quad  &Y_0\in C_b(\mathbb{S}) \quad \stackrel{{\cal L}}{\longleftarrow} \quad \{U_t :0\leq
		t\leq T\} \in L^2_\clZ\big([0,T];\Re^m\big)\quad  &\text{(control input)} &&
	\end{flalign*}
\hrule
\hrule
\end{table*}

\subsection{Stabilizability and detectability}\label{sec:stab-detec}
%

Consider
the solution $\{\mu_t\in \clY^\dagger: t\ge 0\}$ to the forward Kolmogorov
equation:
\begin{equation*}\label{eq:FKE}
	\mu_t(f) = \mu_0(f) + \int_0^t \mu_s\big(\clA f\big)\ud s,
        \qquad t\geq 0
\end{equation*}
The {\em stable complement} of $\clA$ is defined as follows:
\[
S_s := \big\{\mu_0\in \clY^\dagger: \mu_T(f)\,\to\,0 \text{ as }T\to \infty, \forall\, f\in \clY\big\}
\]
Observe that a constant function is $\clA$-invariant and therefore 
$\mu_T(\ones) = \mu_0(\ones)$. Consequently, $S_s\subset \clM_0(\bS)$. 
It is natural to define stabilizability and detectability as dual properties as follows:
%

\medskip

\begin{definition}\label{def:stabilizability}
	The BSDE~\eqref{eq:dual-bsde} is \emph{stabilizable} if $\Rsp(\clL)^\bot \subset S_s$. 
\end{definition}

\medskip

\begin{definition}\label{def:detectability}
	The HMM $(\clA, h)$ is \emph{detectable} if $\Nsp(\clL^\dagger) \subset S_s$. 
\end{definition}

\medskip

\begin{corollary}\label{cor:detect-stab}
	The HMM $(\clA, h)$ is detectable if and only if the 
	BSDE~\eqref{eq:dual-bsde} is stabilizable.
\end{corollary}

\medskip

\begin{remark}\label{rm:vh-detectabilty}
	We compare with the detectability definition
        in~\cite[Definition V.1]{van2010nonlinear}. An HMM is said to
        be detectable if for any $\mu,\nu\in\clP(\bS)$, 
	\[
	\text{Either} \;\;\sP^\mu|_{\clZ_T} \neq \sP^\nu|_{\clZ_T} \quad
        \text{or} \;\; \|\mu_T - \nu_T\|_{\text{\tiny TV}} \stackrel{(T\to\infty)}{\longrightarrow}  0
	\]
	By~\Thm{thm:observability-Zakai-relation}, this statement is identical to the Def.~\ref{def:detectability}.
\end{remark}

\medskip

\begin{remark}\label{rm:ergodicity}
	A Markov process is said to be \emph{ergodic} if there exists
        an invariant measure $\bmu$ such that for all $\mu_0 \in
        \clP(\bS)$, $\|\mu_T - \bmu\|_{\text{\tiny TV}} 
        \stackrel{(T\to\infty)}{\longrightarrow}  0$. 
	For an ergodic process, consider $\tilde{\mu}_0 \in 
	\clM_0(\bS)$ and pick $\mu_0^{(1)},\mu_0^{(2)}\in\clP(\bS)$ and 
	$c\in\Re$ such that $\tilde{\mu}_0 = c(\mu_0^{(1)}-\mu_0^{(2)})$. Then
	\[
	\|\tilde{\mu}_T\|_{\text{\tiny TV}}  
	\le c\|\mu_T^{(1)}-\bmu\|_{\text{\tiny TV}}  +c\|\mu_T^{(2)}-\bmu\|_{\text{\tiny TV}}  
	\stackrel{(T\to\infty)}{\longrightarrow} 0
	\]
	Therefore, if the state process is ergodic then $S_s =
        \clM_0(\bS)$ and the BSDE is stabilizable (and from duality the
        HMM is detectable) irrespective of the observation function $h$.
\end{remark}


%
%
%

\subsection{Linear Gaussian model}
\label{rm:duality-reduces-LG}

Consider the linear Gaussian model:
\begin{subequations}\label{eq:linear-Gaussian-model}
	\begin{align}
		\ud X_t &= A^\tp X_t \ud t + \sigma \ud B_t,\quad X_0\sim N(m_0,\Sigma_0) \label{eq:linear-Gaussian-model-a}\\
		\ud Z_t &= H^\tp X_t \ud t + \ud W_t \label{eq:linear-Gaussian-model-b}
	\end{align}
\end{subequations}
where the prior $N(m_0,\Sigma_0)$ is a Gaussian density with mean
$m_0\in \Re^d$ and variance $\Sigma_0 \succeq 0$, and the process
noise $B=\{B_t\in\Re^p:0\leq t\leq T\}$ is a B.M.  It is assumed that
$X_0,B,W$ are mutually independent. The model parameters
$A\in\Re^{d\times d}$, $H\in\Re^{d\times m}$, and
$\sigma\in\Re^{d\times p}$.

We impose the following restrictions:
	\begin{itemize}
		\item The control input $U=u$ is restricted to be a deterministic
		function of time.  In particular, it does not depend upon the
		observations.  
		For such a control input, the solution $Y=y$ of the
                BSDE is a deterministic function of time, and $V=0$.
                The BSDE becomes a PDE:
		\begin{equation}\label{eq:det_pde}
			-\frac{\partial y_t}{\partial t}(x) = (\clA
                        y_t)(x) + x^\tp H 
			u_t,\quad y_T = c\ones 
		\end{equation}
		where the lower-case notation is used to stress the fact that $u$
		and $y$ are now deterministic functions of time.
		
		\item Consider a
                  finite ($d$-) dimensional space of linear functions:
		\[
		{\sf L}:=\{f:\Re^d\to \R \;:\;  f(x)=x^\tp \tilde{f}\;\text{where}\; \tilde{f}\in\Re^d\}
		\]  
		Then ${\sf L}$ is an invariant subspace for the
		dynamics~\eqref{eq:det_pde}.  On ${\sf L}$, expressing
                $y_t(x) = x^\tp \tilde{y}_t$, the PDE reduces to an ODE:  
		\begin{equation}\label{eq:LTI-ctrl}
		-\frac{\ud y_t}{\ud t}
		= A y_t + H u_t,\quad y_T = 0
		\end{equation}
                where we have dropped the tilde for notational ease.
                The terminal condition is set to 0 because it is the
                only constant function which is also linear.
	\end{itemize}
In this manner, we have recovered the dual control system~\eqref{eq:LTI-ctrl-intro} familiar from the linear systems
theory (and discussed in~\Sec{sec:introduction}).  The main assumption
in going from BSDE~\eqref{eq:dual-bsde} to the ODE~\eqref{eq:LTI-ctrl}
is that the control $U$ is deterministic.  Evidently, such a
choice suffices for the purposes of linear Gaussian estimation. 
A detailed explanation of why this is the case appears in part II.  

\medskip

\begin{remark}
In the stability proofs of Kalman filter, detectability of the
deterministic linear system~\eqref{eq:LTI-obs-intro} is a standard
condition used in proofs (the condition is both necessary and sufficient to show that the
solution of the differential Riccati equation converges~\cite[Thm.~3.7]{kwakernaak1969linear}).
Because the deterministic and stochastic models are different
(compare~\eqref{eq:LTI-obs-intro}
and~\eqref{eq:linear-Gaussian-model}), the fact that detectability of
the deterministic model works also for the stochastic model is somewhat surprising.  The above 
provides an explanation by showing that 
the dual control systems for the two models are the same
(compare~\eqref{eq:LTI-ctrl-intro} and~\eqref{eq:LTI-ctrl}).  
\end{remark}


\subsection{Comparison between linear and nonlinear systems}

For pedagogical reasons, it is useful to draw parallels between the
linear deterministic and the nonlinear stochastic cases.  In both
cases, controllability (resp. observability) is a property of the
range space (resp. null space) of a certain linear operator
(resp. the adjoint operator).  Table~\ref{tab:tab1} provides a comparison of the
linear operators together with the domain and the co-domain spaces.
Based on these definitions, Table~\ref{tb:comparison} provides a
side-by-side comparison of the controllability-observability duality
in the two cases.        

\begin{table*}[t]
	\centering
	\renewcommand{\arraystretch}{1.7}
	\caption{Comparison of the controllability--observability duality for linear and nonlinear systems} \label{tb:comparison}
	\small
	\begin{tabular}{p{0.14\textwidth}p{0.37\textwidth}p{0.4\textwidth}}
		& {\bf Linear deterministic systems} & {\bf Nonlinear stochastic systems} \\ \hline \hline
		Function space for inputs and outputs & $\clU =  L^2([0,T];\Re^m)$\vspace{4pt} \newline  $\langle u,v\rangle = \displaystyle\int_0^T u_t^\tp v_t\ud t$ \vspace{3pt} & $\clU =  L^2_\clZ(\Omega\times [0,T];\Re^m)
		$\vspace{4pt} \newline $\langle U,V\rangle = \displaystyle\hE\Big(\int_0^T U_t^\tp V_t\ud t\Big)$\vspace{3pt} \\ \hline
		Function space for the dual state & $\clY =\Re^d$ \vspace{3pt} \newline $\langle x,y \rangle = x^\tp y$ & $\clY = C_b(\bS)$, $\clY^\dagger = {\cal M}(\bS)$ \vspace{3pt} \newline $\langle \mu,y\rangle = \mu(y)$ \\ \hline
		Linear operators & $\clL : \clU \to \clY$, $u\mapsto y_0$ by ODE~\eqref{eq:LTI-ctrl-intro} \vspace{2pt} \newline 
		$\clL^\dagger : \clY \to \clU$, 
                  \vspace{2pt} \newline \phantom{$\clL^\dagger$}
                  \hspace{0.15em} $x_0\mapsto
                                  \{z_t:0\leq t\leq T\}$ by ODE~\eqref{eq:LTI-obs-intro}
		& $\clL:\clU\times\Re\to\clY$, $(U,c)\mapsto Y_0$ by BSDE~\eqref{eq:dual-bsde} \vspace{2pt}\newline
		$\clL^\dagger:\clY^\dagger\to \clU\times\Re$,
                  \vspace{2pt} \newline \phantom{$\clL^\dagger$}
                  \hspace{0.15em} $\mu\mapsto
                  (\{\sigma_t(h):0\leq t\leq T\}, \mu(\ones))$ by Zakai Eq.~\eqref{eq:Zakai}
		\\ \hline
		Controllability & $\Rsp(\clL) = \Re^d$ & $\overline{\Rsp(\clL)} = \clY$ \\ \hline
		Observability & $\Nsp(\clL^\dagger) = \{0\}$ & $\Nsp(\clL^\dagger) = \{0\}$ \\ \hline
		Duality & \multicolumn{2}{l}{$\Rsp(\clL)^\bot =
                          \Nsp(\clL^\dagger) \quad \Longrightarrow
                          \quad$ A state-output system is observable iff
                          the dual control system is controllable.
	}
		\\ \hline \hline
	\end{tabular}
\end{table*}

\section{Explicit formulae for finite state-space} \label{sec:finite-observability}

\subsection{Notation}

The state-space $\bS$ is finite, namely $\bS= \{1,2,\ldots,d\}$.  In
this case, the space $C_b(\bS)$ and $\clM(\bS)$ are both isomorphic to
$\Re^d$: a real-valued function $f$ (resp. a measure $\mu$) are both
identified with a column vector in $\Re^d$ where the $i^{\text{th}}$
element of the vector represents $f(i)$ (resp. $\mu(i)$), and $\mu(f)
= \mu^\tp f$.  In this manner, the
observation function $h$ is also identified with a matrix
$H\in\Re^{d\times m}$.  We denote the $j^{\text{th}}$ column and the
$i^{\text{th}}$ row of the matrix $H$ by $H^j$ and $H_i$, respectively. 

The generator $\clA$ of the Markov process is identified with a
row-stochastic rate matrix $A\in\Re^{d\times d}$ (the non-diagonal
elements of $A$ are non-negative and the row sum is zero). It acts on
a function through: $
\clA: f\mapsto A f$. 
Its adjoint $\clA^\dagger$ acts on measures through: 
$\clA^\dagger:\mu \mapsto A^\tp \mu$ where $A^\tp$ is the matrix
transpose. 


\subsection{Dual control system}

The dual processes $Y$ and
$V$ are $\Re^d$ and $\Re^{d\times m}$-valued, respectively.  The BSDE~\eqref{eq:dual-bsde} is
finite-dimensional as follows:
\begin{align}
	-\ud Y_t &= \Big(AY_t +HU_t+ \sum_{j=1}^m\dv(H^j) V_t^j\Big)\ud
  t - V_t \ud Z_t,
\nonumber \\
\label{eq:dual-ctrl-finite}
\quad Y_T &= c\ones
\end{align}
where $\ones$ is now the $d$-dimensional column vector of all ones, $\dv(H^j)$ is
a diagonal matrix formed from the $j^{\text{th}}$ column of the matrix $H$, 
and $V_t^j$ denotes the $j^{\text{th}}$ column of the matrix $V_t$.
The solution pair is $(Y,V) \in L^2_\clZ([0,T];\Re^d)\times
L^2_\clZ([0,T];\Re^{d\times m})$.

We refer to the BSDE~\eqref{eq:dual-ctrl-finite} as
the nonlinear model $(A,H)$ and the ODE~\eqref{eq:LTI-ctrl-intro}
as the linear model $(A,H)$.

\subsection{Controllable subspace}

The controllable space $\clC$ is a subspace of $\Re^d$.  Note that if
$U$ is deterministic then $V=0$ and the
BSDE~\eqref{eq:dual-ctrl-finite} reduces to the
ODE~\eqref{eq:LTI-ctrl-intro}.  This means that the controllable
subspace for the linear model $\sp\{H, \,  AH, \, A^2H, \ldots \}
\subset \clC$.  Directly by
applying Prop.~\ref{thm:controllable-subspace}, 
\begin{equation}\label{eq:obs_gram_nl}
\begin{aligned}
	\clC = \sp\big\{\ones, &\,  H, \,  AH, \,  A^2H, \,  A^3H, \, \ldots, \\
	&H\cdot H, \,  A(H \cdot H), \,  H\cdot (AH), \,  A^2(H\cdot H),\ldots, \\
	&H\cdot (H\cdot H), \,  (AH)\cdot (H\cdot H), \, \ldots \big\} 
\end{aligned}
\end{equation}
where the dot notation denotes the element-wise (Hadamard) product
between two matrices.  Again, one observes that the first row $\{H,
AH, A^2H,\ldots \}$ is the same as the controllability matrix of the
linear model $(A,H)$.
Therefore, if $(A,H)$ is controllable in the sense of linear systems
theory then the nonlinear model is also controllable (from duality 
the HMM is then observable.)  The presence of additional entries
in~\eqref{eq:obs_gram_nl} means
that the nonlinear model $(A,H)$ may be controllable  even though the
linear model $(A,H)$ is not.  
To highlight the difference, the following proposition gives a
sufficient condition for controllability which does not depend upon $A$. 
\medskip

\begin{proposition}\label{prop:sufficient}
	The nonlinear model $(A,H)$ is controllable if $H$ is an injective map from
	$\bS$ into $\Re^m$ (the map is injective if and only if $H_i \neq H_j$ for all $i\neq
	j$).    
	If $A=0$ then the injective property is
	also necessary for controllability.    
\end{proposition}
\medskip
\begin{proof}
See \Appendix{ss:pf-prop45}.
\end{proof}
\medskip
\begin{remark}\label{rm:finite-case-equivalency}
In~\cite{van2009observability}, a test for stochastic observability is
given for the white noise observation model in the finite state-space
settings.  The test is given in terms of the dimension of the space of
observable functions
(Def.~\ref{def:un-observable-measures-observable-functions}). Since
the notation is somewhat different, we recall
that~\cite{van2009observability} denotes the set of distinct possible
values of noise-free observations by $h(\bS) := \{h_1,\ldots,h_r\}$
(the value $h_i$ should not to be confused with $h(i)$).    
Note $r \le d$ with
equality holds when $h(i)\neq h(j)$ for all $i\neq j$.  
Next for each $h_k\in h(\bS)$, a diagonal projection matrix $P_{h_k}
\in\Re^{d\times d}$ is defined whose non-zero elements are
$P_{h_k}(i,j) = 1$ if $h(i)=h_k$ and $i=j$.  In terms of these matrices, the space of
observable functions is shown to be~\cite[Lemma 9]{van2009observability}
	\[
	\clO = \sp\big\{P_{n_0}AP_{n_1}AP_{n_2}\cdots AP_{n_k}\ones:k\ge 0,n_i\in h(\bS)\big\}
	\]
	It is shown in \Appendix{ssec:justifiction-vanHandel-and-me-finite}
        that $\clO= \clC$ (formula in~\eqref{eq:obs_gram_nl}).
	%
	
\end{remark}	

\subsection{Controllability gramian}

For the nonlinear model $(A,H)$, the controllability gramian $\sW$ 
is most directly expressed in
terms of the solution operator of the Zakai equation defined as follows:
\[
\ud \Psi_t = A^\tp \Psi_t \ud t + \sum_{j=1}^m\dv(H^j)\Psi_t\ud Z_t^j,\quad \Psi_0 = I_d
\]
where $I_d$ is the $d\times d$ identity matrix.  
In \Appendix{ssec:derivation-ctrl-gramian}, it is shown that
\begin{equation}\label{eq:gramian-finite}
\sW = \ones\ones^\tp + \tE\Big(\int_0^T \Psi_t^\tp HH^\tp \Psi_t \ud t\Big)
\end{equation}
Since $\sW$ is a deterministic matrix, controllability
admits to a rank condition test:
\[
(A,H) \; \text{is controllable} \quad \Longleftrightarrow \quad \sW \text{ is full rank} 
\]

\medskip

\subsection{Stabilizability}

Because $A$ is a stochastic matrix, a simple application of the Ger\v{s}gorin
circle theorem shows that the eigenvalues of $A$ are either in the
open left-half-plane or at zero, and 
%
\[
S_s^\bot = S_0 := \{f \in \Re^d : \; Af = 0\} 
\]
Therefore, stabilizability is equivalent to an inclusion property:
\[
(A,H) \; \text{is stabilizable} \quad \Longleftrightarrow \quad S_0\subset\clC
\]







A meaningful characterization is possible by partitioning the finite state-space $\bS$ into $r$ \emph{ergodic classes}
$\{\bS_k:k=1,2,\ldots,r\}$ as follows: 
\begin{enumerate}
	\item $\bS = 
	\bigcup_{k=1}^r \bS_k$ where $\sP([X_t\in \bS_l] \mid [X_0\in \bS_k]) = 0$ for all $t\ge 0$ 
	and $l \neq k$.
	\item By choosing an appropriate coordinate, the matrix
	\[
	A = \begin{pmatrix}
		A_1 & 0 & \cdots & 0\\
		0 & A_2 & \cdots & 0\\
		\vdots & \vdots & \ddots & \vdots\\
		0 & 0 & \cdots & A_r
	\end{pmatrix}
	\]
	where $A_k$ is a row stochastic matrix on $\bS_k$ for $k = 1,2,\ldots,r$.
\end{enumerate} 
The Markov process is said to be \emph{ergodic} if it has a single
ergodic class ($r = 1$). 
The following proposition provides an explicit characterization of
stabilizability:

\medskip

\begin{proposition}\label{prop:detectability-implications}
	Consider the nonlinear model $(A,H)$ and an associated
        ergodic partition $\bS=\cup_{k=1}^r \bS_k$. Then
	\begin{itemize}
		\item If $\bS$ has a single ergodic class ($r=1$) then $(A,H)$ is stabilizable.  
		\item If $r>1$ then $(A,H)$  is stabilizable 
		if and only if the indicator functions $\ones_{\mathbb{S}_k} \in \clC$ for
		$k=1,2,\ldots,r$.  
	\end{itemize}
\end{proposition}

\medskip

\begin{proof}
The matrix $A$ has exactly $r$ eigenvalues at zero with an $r$-dimensional eigenspace $S_0 =
\sp\{\ones_{\mathbb{S}_k}:k=1,2,\ldots,r\}$.   
\end{proof}

\medskip

\begin{remark}
The two bullets in Prop.~\ref{prop:detectability-implications}
correspond to the ergodic signal case and the non-ergodic signal case
as discussed in \Sec{sec:introduction}.  The first bullet says that if
the Markov process is ergodic then the dual control system is
stabilizable irrespective of $H$.  For the non-ergodic case, the
second bullet provides a simple condition for stabilizability.  It can be shown
that the condition is both necessary and sufficient for the optimal
filter to asymptotically detect the correct ergodic
class~\cite[Ch.~8]{JinPhDthesis} (see also~\cite{baxendale2004asymptotic}).
\end{remark}

\section{Conclusions and directions for future work}
\label{sec:conc}

In the study of deterministic linear and nonlinear systems, duality
has played a central role in defining, interpreting, and using the
property of observability.  In this paper, we presented the first such dual control
system for studying stochastic observability of HMMs with white noise
observations.  We related the dual control system to both nonlinear
filtering (the Zakai equation) and the linear Gaussian model.
The latter relationship is shown to recover the classical duality between
controllability and observability of linear systems.  In fact, the development is
entirely parallel in the nonlinear and the linear systems.  This is
shown with the aid of the Table~\ref{tb:comparison} with a
side-by-side comparison.  


Because the stress in this paper was on the duality between controllability and observability, we did
not explicitly relate the process $Y$ and the hidden process $X$.  In
fact, it is possible as well and yields the following useful formula
\[
\E \big(Y_T(X_T)\big) = \mu(Y_0) - \E \Big( \int_0^T U_t^\tp \ud Z_t\Big)
\]  
The formula is the starting point for the companion paper (part II)
which is concerned with the use of duality to express the nonlinear
filtering problem as a dual optimal control problem.  The connection
to optimal control theory opens up several avenues of research, namely
study of asymptotic stability of the optimal filter, a definition of
suitable supply rates that yields useful dissipative characterizations
of the HMM, and opportunities for algorithm design.  These are also
discussed at length as part of the conclusions of the companion paper.


A natural question also is to relate this work to the study
of deterministic models.  Unfortunately, our work crucially relies on
additive Gaussian form of the measurement noise.  This effectively
precludes the zero measurement noise case.  On the other hand,
considering the limit as the covariance of the measurement noise goes to zero
may be meaningful from robustness perspective, and may yield
useful insights for the deterministic models.  

\section{Acknowledgement}

It is a pleasure to acknowledge Sean Meyn and Amirhossein Taghvaei for many useful
technical discussions over the years on the topic of duality.  The authors also acknowledge Alain Bensoussan for his early encouragement of this work.

\appendix

\section{Proofs of the statements}

%

\subsection{Proof of~\Thm{thm:observability-Zakai-relation}}\label{pf-obs-Zakai}

\begin{itemize}
\item{(1 $\Longleftrightarrow$ 2)} It directly follows from the
relative entropy formula~\eqref{eq:clark-result}.

\item{(2 $\Longrightarrow$ 3)} The Zakai equation~\eqref{eq:Zakai}
with $f=\ones$ gives
\begin{equation}\label{eq:zakai_for_1}
\sigma_t(\ones) = 1 + \int_0^t \sigma_s(h)^\tp\ud Z_s = 1 + \int_0^t \sigma_s(\ones)\pi_s(h)^\tp\ud Z_s
\end{equation}
where the formula~\eqref{eq:normalize-Zakai} is used to obtain
the second equality. It follows that if $\pi_s^\mu(h) = \pi_s^\nu(h)$ for all $0\le s\le t$ then $\sigma_t^\mu(\ones) = \sigma_t^\nu(\ones)$. Using~\eqref{eq:normalize-Zakai}, $\sigma_t^\mu(h) = \sigma_t^\nu(h)$.

\item{(3 $\Longrightarrow$ 2)} From the first equality
in~\eqref{eq:zakai_for_1}, it follows that if $\sigma_s^\mu(h) =
\sigma_s^\nu(h)$ for all  $0\le s\le t$ then $\sigma_t^\mu(\ones) = \sigma_t^\nu(\ones)$. Using~\eqref{eq:normalize-Zakai}, $\pi_t^\mu(h) = \pi_t^\nu(h)$. 
\end{itemize}

\subsection{Proof of~\Thm{thm:observability-definition}}\label{ssec:pf-observability-definition}

The space $\clU=L^2_{\clZ}([0,T];\Re^m)$ is a Hilbert space with the inner product
$\langle U,V\rangle_{\clU} := \tE\Big(\int_0^T U_t^\tp V_t \ud
t\Big)$.  Therefore, $\clU\times \Re$ is a Hilbert space with a
natural inner product $
\langle (U,c),(V,d)\rangle_{\clU\times \Re} := \langle U,V\rangle_{\clU} + cd$.  For a function $f\in\clY$ and a measure $\mu \in \clY^\dagger$,
the duality pairing is denoted by $\langle \mu, f\rangle_\clY :=
\mu(f)$.  Let $U\in \clU$, $c\in \Re$, and $\mu \in \clY^\dagger$. By
linearity, $\clL(U,c) = \clL(U,0)+c\ones$ and therefore
\[
\langle \mu,\clL(U,c) \rangle_\clY  =  \langle \mu,\clL(U,0) \rangle_\clY  + c \mu(\ones)
\]
Thus, the main calculation is to show $\langle
\mu,\clL(U,0)\rangle_\clY = \langle \sigma(h),U\rangle_{\clU}$ where
$\sigma(h)=\{\sigma_t(h)\in\Re^m:0\leq t\leq T\}$ solves the Zakai
equation~\eqref{eq:Zakai} with $\sigma_0=\mu$. 
This is done by using the It\^o-Wentzell formula for
measure valued processes~\cite[Thm.~1.1]{krylov2011ito} (note here
that all stochastic processes
are forward adapted),
\begin{align*}
	\ud \big(\sigma_t(Y_t)\big) &= \big(\sigma_t(\clA Y_t) \ud t + \sigma_t(h^\tp Y_t)\ud Z_t\big) + \sigma_t(h^\tp V_t)\ud t \\
	&\quad+ \big(\sigma_t(-\clA Y_t - h^\tp U_t - h^\tp V_t) \ud t + \sigma_t(V_t)\ud Z_t\big)\\
	&= -U_t^\tp \sigma_t(h)\ud t + \sigma_t(h^\tp Y_t + V_t^\tp) \ud Z_t
\end{align*}
Integrating both sides,
\[
\sigma_T(Y_T) - \mu(Y_0) = -\int_0^T U_t^\tp\sigma_t(h)\ud t + \int_0^T \sigma_t(h^\tp Y_t+V_t^\tp) \ud Z_t
\]
With $Y_T=0$, because $Z$ is a $\tsP$-B.M.,
\begin{equation*}\label{eq:pf_thm1_1}
\langle \mu,\clL(U,0)\rangle_\clY = \mu(Y_0) = \tE\Big(\int_0^T U_t^\tp \sigma_t(h)\ud t \Big) =
        \langle \sigma(h),U\rangle_{\clU}
\end{equation*}
In summary, 
\[
\langle \mu,\clL(U,c) \rangle_\clY = \langle \sigma(h), U \rangle_\clU + c\mu(\ones) = \big\langle \clL^\dagger \mu, (U,c)\big\rangle_{\clU\times \Re}
\]

\subsection{Proof of~\Prop{thm:controllable-subspace}}\label{ss:pf-prop43}

The proof is adapted from~\cite[Thm.~3.2]{peng1994backward}.  For notational ease, we assume $m=1$.  
The definition of $\Nsp(\clL^\dagger)$ is:
\[
\mu \in \Nsp(\clL^\dagger) \Leftrightarrow \mu(\ones) = 0 \text{
	and } \sigma_t(h) \equiv 0 \quad \forall\;t\in[0,T]
\]
Since $\Nsp(\clL^\dagger)$ is the annihilator of $\Rsp(\clL)$, we have
$\ones,h \in \Rsp(\clL)$. Consider next the 
Zakai equation~\eqref{eq:Zakai} with the initial condition
$\mu\in\Nsp(\clL^\dagger)$ and $f=h$:
\[
\sigma_t(h) = \mu(h) + \int_0^t \sigma_s(\clA h) \ud s + \int_0^t \sigma_s(h^2) \ud Z_s
\]
Since $t$ is arbitrary, the left-hand side is identically zero for all $t\in[0,T]$ if and only if
\[
\mu(h) = 0,\quad \sigma_t(\clA h) \equiv 0,\quad \sigma_t(h^2)\equiv 0 \quad \forall\;t\in[0,T]
\]
and in particular, this implies $\clA h, h^2\in\Rsp(\clL)$. 

The subspace $\clC$ is obtained by continuing to repeat the steps
ad infinitum: If at the conclusion of the  $k^\text{th}$ step, we find
a function $g\in \clC$ such that $\sigma_t(g)\equiv 0$ for all
$t\in[0,T]$.  Then through the use of the Zakai equation,
\[
\mu(g) = 0,\quad \sigma_t(\clA g) \equiv 0,\quad \sigma_t(hg)\equiv 0
\quad \forall\;t\in[0,T]
\] 
so $\clA g, hg \in \clC$.  By construction, because $
\mu \in \Nsp(\clL^\dagger)$, $\clC
= \Rsp(\clL)$. 

\subsection{Proof of~\Prop{prop:gramian-observability}}\label{ss:pf-prop44}


Suppose $f \in \Rsp(\sW)$. Then there exists $\mu\in\clM(\bS)$ such that $\sW\mu = f$. Let $\big(U,\mu(\ones)\big) = \clL^\dagger\mu$, and apply the control $U$ to the BSDE with terminal condition $Y_T = \mu(\ones)\ones$. Then 
\[
Y_0 = \clL(U,\mu(\ones)) = \clL \clL^\dagger \mu = \sW\mu = f
\]
Suppose another $(\tilde{U},c)\in\clU\times \Re$ gives $\clL(\tilde{U},c) = f$. Then
\begin{align*}
0 &= \big\langle \mu,\clL\big(U-\tilde{U},\mu(\ones)-c\big) \big\rangle_\clY\\
&= \big\langle \clL^\dagger\mu , \big(U-\tilde{U},\mu(\ones)-c\big)
  \big\rangle_{\clU\times \Re} \\
&= \big\langle \big(U,\mu(\ones)\big), \big(U-\tilde{U},\mu(\ones)-c\big) \big\rangle_{\clU\times \Re} 
\end{align*}
The minimum norm formula~\eqref{eq:min-norm} follows because of this
orthogonality property. 

\subsection{Proof of Proposition~\ref{prop:sufficient}}\label{ss:pf-prop45}

\newP{Step 1} We first provide the proof for the case when $m=1$. In
this case, $H$ is a column vector and $H_i$ denotes its
$i^{\text{th}}$ element. 
We claim that if $H_i\neq H_j$ for all $i\neq j$, then 
\begin{equation}\label{eq:subset}
	\sp\{\ones, \, H, \, H\cdot H, \, \ldots,\, \underbrace{H\cdot H \cdots H}_{(d-1)\text{ times}}\} = \Re^d
\end{equation}
where (as before) the dot denotes the element-wise product.
Assuming that the claim is true, the result easily follows because the
vectors on left-hand side are contained in $\Rsp(\clL)$
(see~\eqref{eq:obs_gram_nl}). It remains to prove the claim. 
For this purpose, express the left-hand side of~\eqref{eq:subset} as
the column space of the following matrix:
\begin{align*}
	\begin{pmatrix}
		1 & H_1 & H_1^2 & \cdots & H_1^{d-1}\\
		1 & H_2 & H_2^2 &\cdots & H_2^{d-1}\\
		\vdots &\vdots  & \vdots &\cdots & \vdots\\
		1 & H_d & H_d^2 & \cdots & H_d^{d-1}
	\end{pmatrix}
\end{align*}
This matrix is easily seen to be full rank by using the Gaussian elimination:
\[
\begin{pmatrix}
	1 & H_1 & H_1^2 & \cdots & H_1^{d-1}\\
	0 & H_2-H_1 & H_2^2-H_1^2 &\cdots & H_2^{d-1}-H_1^{d-1}\\
	\vdots &\vdots  & \vdots &\cdots & \vdots\\
	0 & 0 & 0 & \cdots & \prod_{i=1}^{d-1} (H_d-H_i)\\
\end{pmatrix}
\]
The diagonal elements are non-zero because $H_i\neq H_j$.

\newP{Step 2} In the general case, $H$ is a $d \times m$ matrix and
$H_i$ denotes its $i^\text{th}$ row.  We claim that if $H_i\neq H_j$
for all $i\neq j$ then there exists a vector $\tilde{H}$ in the column
span of $H$ such that $\tilde{H}_i\neq \tilde{H}_j$ for all $i\neq
j$. Assuming that the claim is true, the result follows from the $m=1$
case by considering~\eqref{eq:subset} with $\tilde{H}$.
It remains to prove the claim.  Let $\{e_1,\ldots,e_d\}$ denote the
canonical basis in $\Re^d$.  The assumption means $(e_i-e_j)^\tp H$
is a non-zero row-vector in $\Re^m$ for all $i\neq j$. 
Therefore, the null-space of $(e_i-e_j)^\tp H$ is a
$(m-1)$-dimensional hyperplane in $\Re^m$. Since there are only  
finite such hyperplanes, there must exist a vector $a\in\Re^m$
such that $(e_i-e_j)^\tp Ha \neq 0$ for all $i\neq j$. Pick such an
$a$ and define $\tilde{H} := Ha$.

\newP{Step 3} To show the necessity of the injective property when
$A=0$, assume $H_i=H_j$ for some $i\neq j$. Then the corresponding rows
are identical, so the rank is less than $d$.


\subsection{Proof that $\clO=\clC$ in Rem.~\ref{rm:finite-case-equivalency}} \label{ssec:justifiction-vanHandel-and-me-finite}

We begin with $A=0$ case so
\[
\clC = \sp\{\ones, H, \dv(H)H, \dv(H)^2H,\ldots\}
\]
Since $\dv(H)^nH = [h^{n+1}(1),\ldots,h^{n+1}(d)]^\tp$, an element of $f\in \clC$ can be expressed by
\[
f=\sum_{j=0}^\infty a_j[h^j(1),\ldots,h^j(d)]^\tp =  \sum_{k=1}^r\sum_{j=0}^\infty a_j h_k^j P_{h_k}\ones \in \clO
\]
Therefore, $\clC \subset \clO$. To show $\clO\subset\clC$, let
\[
f=\sum_{k=1}^r b_kP_{h_k}\ones
\]
It suffices to show that there exists $\{a_j:j=0,1,\ldots\}$ such that
$b_k = \sum_{j=0}^\infty a_j h_k^j$ for all $k=1,\ldots,r$. In fact,
such $a_j$ can be found explicitly by setting $a_j = 0$ for $j\geq r$ and inverting the following matrix:
\[
\begin{pmatrix}
	1 & h_1 & h_1^2 &\cdots & h_1^{r-1}\\
	1 & h_2 & h_2^2 &\cdots & h_2^{r-1}\\
	\vdots &\vdots & \vdots &\cdots & \vdots\\
	1 & h_r & h_r^2 &\cdots & h_r^{r-1}\\
\end{pmatrix}
\]
It is invertible because $h_i$ are distinct (the proof is by Gaussian
elimination as before).


\medskip

For general $A\neq 0$ case, we repeat the same procedure as above for arbitrary matrices $M_1$ and $M_2$ which are multiples of $A$ and $\dv(H)$, to claim that
\begin{align*}
	\sp\{M_1M_2H, M_1&\dv(H)M_2H, M_1\dv(H)^2M_2H,\ldots\} \\
	&= \sp\{M_1P_{h_k}M_2H:k=1,\ldots,r\}
\end{align*}
The proposition is proved by repeating this for countable times.

\subsection{Formula~\eqref{eq:gramian-finite} for the gramian} \label{ssec:derivation-ctrl-gramian}

For a given measure $\mu\in\Re^d$, $\sW\mu=Y_0$ is obtained by solving
the BSDE
\begin{align*}
-\ud Y_t &= \big(AY_t + HH^\tp \sigma_t + \sum_{j=1}^m H^j\cdot
  V_t^j\big) \ud t - V_t \ud Z_t \\
Y_T & = (\mu^\tp\ones)\ones
\end{align*}
Consider the process
\[
\Theta_t:= \Psi_t^\tp Y_t + \int_0^t \Psi_s^\tp HH^\tp \sigma_s \ud s,\quad 0\le t \le T
\]
Then by the It\^o product formula,
\[
\ud \Theta_t = \sum_{j=1}^m\Psi_t^\tp \big(H^jY_t + V_t^j\big) \ud Z_t^j
\]
Therefore, $\{\Theta_t:0\leq t\leq T\}$ is a $\tsP$-martingale. In particular, 
\[
Y_0 = \Theta_0 = \tE(\Theta_T)=\tE\big(\Psi_T^\tp 11^\tp \mu + \int_0^T \Psi_t^\tp HH^\tp\sigma_t\ud t\big) 
\]
Since the un-normalized filter is given by $\sigma_t = \Psi_t\mu$, 
\[
\sW \mu = \tE\Big(\Psi_T^\tp 11^\tp + \int_0^T \Psi_t^\tp H H^\tp\Psi_t\ud t \Big) \mu
\]
Finally, $\tE(\Psi_T^\tp \ones\ones^\tp) = \ones\ones^\tp$ because
$\frac{\ud}{\ud t} \tE(\Psi_t^\tp \ones) = 0$.

\bibliographystyle{IEEEtran}
\bibliography{../bibfiles/_master_bib_jin.bib,../bibfiles/jin_papers.bib,../bibfiles/extrabib.bib,../bibfiles/estimator_controller.bib}

 \begin{IEEEbiography}[{\includegraphics[width=1in,height=1.25in,clip,keepaspectratio]{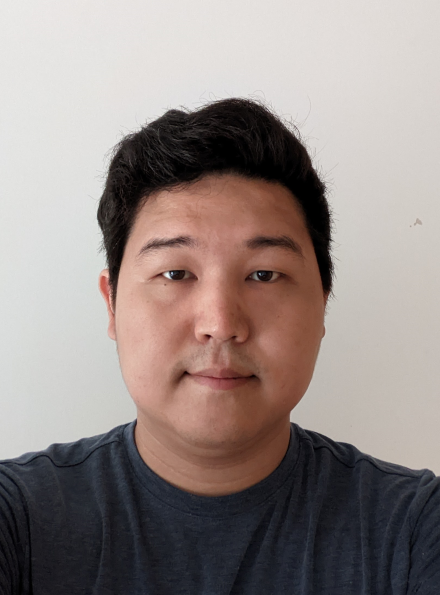}}]{Jin Won Kim} received the Ph.D. degree in Mechanical Engineering from University of Illinois at Urbana-Champaign, Urbana, IL, in 2022.
 	He is now a postdocdoral research scientist in the Institute of Mathematics at the University of Potsdam.
 	His current research interests are in nonlinear filtering and stochastic optimal control.
 	He received the Best Student Paper Awards at the IEEE Conference on Decision and Control 2019.
 \end{IEEEbiography}

\begin{IEEEbiography}[{\includegraphics[width=1in,height=1.25in,clip,keepaspectratio]{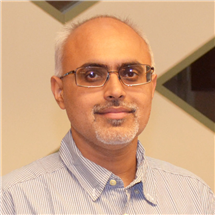}}]{Prashant G. Mehta} received the Ph.D. degree in Applied Mathematics from Cornell University, Ithaca, NY, in 2004.
	He is a Professor of Mechanical Science and Engineering at the University of Illinois at Urbana-Champaign.
	Prior to joining Illinois, he was a Research Engineer at the United Technologies Research Center (UTRC). His current research interests are in nonlinear filtering. He received the Outstanding Achievement Award at UTRC for his contributions to the modeling and control of combustion instabilities in jet-engines. His students received the Best Student Paper Awards at the IEEE Conference on Decision and Control 2007, 2009 and 2019, and were finalists for these awards in 2010 and 2012. In the past, he has served on the editorial boards of the ASME Journal of Dynamic Systems, Measurement, and Control and the Systems and Control Letters. He currently serves on the editorial board of the IEEE Transactions on Automatic Control. 	
 	\end{IEEEbiography}

\end{document}